\documentclass[francais]{article} 
\usepackage{amsmath,amsthm,amsfonts,amssymb,amscd}
\usepackage{babel}
\usepackage{graphicx}
\evensidemargin = \oddsidemargin
\advance\textwidth by -2in
\theoremstyle{plain}
\newtheorem{theoreme}{Th\'eor\`eme}
\newtheorem{proposition}[theoreme]{Proposition} 
\newtheorem{lemme}[theoreme]{Lemme}
\newtheorem{corollaire}[theoreme]{Corollaire} 

\theoremstyle{definition}
\newtheorem{definition}[theoreme]{D\'efinition} 
\theoremstyle{remark}
\newtheorem*{remarque}{Remarque}

\newtheorem*{remerciements}{Remerciements}

\def\Alinea#1{\hfill\break%
  \hbox to \parindent{\hss{{\rm #1}}\enspace}\ignorespaces}
\def\alinea#1{\par\noindent%
  \hbox to \parindent{\hss{{\rm #1}}\enspace}\ignorespaces}
\gdef\up{\textup}
\def\from{\colon}
\def\res#1{\,\vert\,{}_{#1}}
\def\classe#1{\mathcal{C}_{}^{#1}}
\def\Chi{\setbox0=\hbox{$\chi$} \mathord{\raise\dp0\hbox{$\chi$}}}
\def\eps{\varepsilon}

\def\D{\mathbf{D}}

\def\R{\mathbf{R}}
\def\S{\mathbf{S}}
\def\T{\mathbf{T}}
\def\Z{\mathbf{Z}}
\def\CP^#1{\mathbf{P}^{#1}(\mathbf{C})}
\def\RP^#1{\mathbf{P}^{#1}(\mathbf{R})}
\def\id{\mathrm{id}}
\def\CC{\mathcal{C}}
\def\DD{\mathcal{D}}
\def\FF{\mathcal{F}}
\def\PP{\mathcal{P}}
\def\XX{\mathcal{X}}
\def\SC{\mathcal{SC}}
\def\SCT{\mathcal{SCT}}
\def\O{\mathrm{O}}
\def\SO{\mathrm{SO}}
\def\SL{\mathrm{SL}}

\def\Int{\mathop{{\rm Int}}\nolimits}

\def\op{\mathopen}
\def\cl{\mathclose}
\def\ol{\overline}
\def\Ol#1{\smash{\,\overline{\!#1}}}
\def\wt{\widetilde}
\def\Wt#1{\smash{\,\widetilde{\!#1}}}

\def\adresse#1{\def\ls{\egroup\egroup\hbox\bgroup\itshape\bgroup} \par \noindent
  \hbox to \textwidth {\hfill
  \vbox {\small \hbox \bgroup\itshape\bgroup #1 \egroup\egroup}}}
\title{Sur les transformations de contact\\ au-dessus des surfaces}
\author{Emmanuel \textsc{Giroux}\thanks
{\ Centre National de la Recherche Scientifique (\textsc{umr} 5669).}}
\begin{document}

\maketitle

La g\'eom\'etrie de contact, n\'ee des travaux de C.~Jacobi et de S.~Lie sur les
\'equations diff\'erentielles et leurs sym\'etries, a \'et\'e la g\'eom\'etrie
des transformations de contact avant d'\^etre celle des structures de contact.
On s'int\'eresse ici, dans un cas simple, \`a la topologie du groupe que ces
transformations composent.

Soit $S$ une surface et $\pi \from V \to S$ le fibr\'e des droites coorient\'ees
tangentes \`a~$S$, traditionnellement nomm\'e aussi fibr\'e des \'el\'ements de
contact coorient\'es au-dessus de~$S$. Le champ de plans tautologique~$\xi_S$
d\'efini sur~$V$ par
$$ \xi_S (q,\delta) = \bigl( D\pi(q,\delta) \bigr)^{-1} (\delta), \quad
   \text{o\`u $q \in S$ et $\delta \subset T_qS$,} $$
est l'arch\'etype des structures de contact. Ce champ est en outre
cano\-ni\-quement
coorient\'e --~puisque les \'el\'ements de contact le sont~-- et on appelle {\it
transformation de contact\/} au-dessus de~$S$ tout diff\'eomorphisme de~$V$ qui
pr\'eserve~$\xi_S$ ainsi que sa coorientation. Ces transformations engendrent un
groupe $\DD(V ;\xi_S)$ dans lequel le groupe $\DD(S)$ des diff\'eomorphismes
de~$S$ se plonge naturellement : tout diff\'eomorphisme de~$S$ se rel\`eve \`a
$V$ en une transformation de contact. De m\^eme, le groupe $\DD (V,\partial V ;
\xi_S)$ des transformations de contact relatives au bord --~{\it i.e.} \'egales
\`a l'identit\'e sur $\partial V$~-- contient le groupe $\DD (S,\partial S)$
des diff\'eomorphismes de~$S$ relatifs au bord ou, plus exactement, tangents \`a
l'identit\'e le long du bord. Ces groupes et les autres espaces fonctionnels qui
interviennent plus loin sont munis de la topologie de la convergence uniforme
$\classe\infty$ sur les parties compactes. Le but principal de cet article est
d'\'etablir le r\'esultat suivant :

\begin{theoreme} 
Soit $S$ une surface compacte, connexe, orientable et qui n'est ni une sph\`ere
ni un tore. Les plongements naturels
\begin{align*}
\DD (S, \partial S) & \longrightarrow \DD (V, \partial V ; \xi_S), \\
\DD (S) & \longrightarrow \DD (V ; \xi_S) \\ \text{\llap{et\quad}}
\DD (\Int S) & \longrightarrow \DD (\Int V ; \xi_S)
\end{align*}
sont des $0$-\'equivalences d'homotopie et, plus pr\'ecis\'ement, mettent en
bijection les composantes connexes de la source et du but.
\end{theoreme}

En fait, les plongements ci-dessus sont tr\`es probablement des \'equivalences
d'homotopie. Au prix de complications surtout techniques, la m\'e\-thode suivie
ci-apr\`es semble d'ailleurs applicable aux familles --~de diff\'eomorphismes,
de plongements, de structures de contact~-- d\'ependant d'un nombre quelconque
de param\`etres.

Dans ce th\'eor\`eme, l'injectivit\'e des applications induites au niveau des
composantes connexes est \'evidente. Des arguments topologiques simples montrent
en effet que les plongements 
\begin{align*}
\DD (S, \partial S) & \longrightarrow \DD (V, \partial V), \\
\DD (S) & \longrightarrow \DD (V) \\ \text{\llap{et\quad}}
\DD (\Int S) & \longrightarrow \DD (\Int V)
\end{align*}
induisent eux-m\^emes des injections. La surjectivit\'e, en revanche, n'\'etait
apparemment pas connue --~par exemple pour le disque~\cite{Ar}~-- et implique le
r\'esultat suivant :

\begin{corollaire} 
Soit $S$ une surface compacte, connexe, orientable et qui n'est ni une sph\`ere
ni un tore. Les plongements naturels
\begin{align*}
\DD (V, \partial V ; \xi_S) & \longrightarrow \DD (V, \partial V), \\
\DD (V ; \xi_S) & \longrightarrow \DD (V) \\ \llap{\text{et\quad}}
\DD (\Int V ; \xi_S) & \longrightarrow \DD (\Int V) 
\end{align*}
appliquent injectivement les composantes connexes de la source dans celles du
but.
\end{corollaire}

Pour compl\'eter le tableau, on d\'ecrit ci-dessous ce qu'il advient lorsque $S$
est une sph\`ere ou un tore.

\begin{theoreme} 
Le groupe des transformations de contact au-dessus de la sph\`ere est
connexe.
\end{theoreme}

\`A titre de test, on remarque que le flot g\'eod\'esique de la sph\`ere entre
les instants $0$ et $\pi$ m\`ene de l'identit\'e \`a la diff\'erentielle de
l'application antipodale par un chemin de transformations de contact.

\medskip

Lorsque $S$ est un tore, la vari\'et\'e de contact $(V, \xi_S)$ s'identifie au
tore $\T^3$ muni de sa structure de contact canonique~$\xi$, dont une \'equation
de Pfaff est par exemple
$$ \cos \theta \, dx + \sin \theta \, dy = 0, \qquad
   (x,y,\theta) \in \R^2 \!/ \Z^2 \times \R / 2\pi \Z \,. $$ 

\begin{theoreme} 
Soit $\Pi$ le sous-groupe de $\SL_3(\Z)$ constitu\'e des transformations qui
laissent invariant $\Z^2 \times \{0\} \subset \Z^3$. Le morphisme
$$ \pi_0 \DD (\T^3,\xi) \longrightarrow \SL_3(\Z) $$
fourni par l'action sur l'homologie est injectif et a pour image~$\Pi$.
\end{theoreme}

\begin{remarque}
Les diff\'eomorphismes de $V$ qui pr\'eservent~$\xi_S$ mais pas sa coorientation
forment, modulo l'action des isotopies de contact, un groupe isomorphe au
produit direct de $\Z/2\Z$ par $\pi_0 \DD (V ;\xi_S)$ --~respectivement par
$\pi_0 \DD (V,\partial V ;\xi_S)$ ou par $\pi_0 \DD (\Int V ;\xi_S)$ si on
consid\`ere les diff\'eomorphismes relatifs au bord ou les diff\'eomorphismes de
l'int\'erieur. Le g\'en\'erateur du facteur $\Z/2\Z$ correspond \`a
l'application antipodale sur les fibres de~$V$.
\end{remarque}

\begin{remerciements}
L'\'etude pr\'esent\'ee ici a pour point de d\'epart des conversations avec
Vladimir~Arnold --~en mars 1994 \`a Paris~-- et avec Ana~Cannas da~Silva --~en
juillet 1997 \`a Park City. Par la suite, lors d'une premi\`ere approche du
probl\`eme bas\'ee sur l'examen des n\oe uds legendriens, j'ai eu plusieurs
discussions utiles avec Maia~Fraser. Je les remercie tous les trois. Je remercie
aussi Leonid~Polterovich pour m'avoir fait part de son travail r\'ecent~\cite{Po}
o\`u le th\'eor\`eme~1 trouve un \'echo.
\end{remerciements}

\section{Plongements et diff\'eomorphismes de contact}

Avant de d\'emontrer les th\'eor\`emes principaux, on fait quelques observations
g\'en\'erales sur la topologie des groupes $\DD (V,\partial V ;\xi)$ et $\DD
(V ;\xi)$ pour une vari\'et\'e de contact $(V,\xi)$ quelconque de dimension~$3$.
Ces remarques s'inscrivent dans le m\^eme esprit que la strat\'egie propos\'ee
par K.~Honda dans \cite{Ho} pour analyser les structures de contact tendues sur les
vari\'et\'es suffisamment grandes. Par convention, les structures de contact
consi\-d\'er\'ees ci-apr\`es sont toujours orient\'ees et les diff\'eomorphismes
de contact pr\'eservent cette orientation.

\medskip

On rappelle que toute surface $F$ plong\'ee dans $(V,\xi)$ h\'erite d'un {\it
feuilletage caract\'eristique}, not\'e $\xi F$, qui est engendr\'e par le champ
de droites singulier $\xi \cap TF$. Dans la suite, $\xi$ et~$F$ sont orient\'ees
de sorte que $\xi F$ l'est aussi. Concr\`etement, si $\xi = \ker \alpha$ et si
$\omega$ est une forme d'aire sur~$F$, le feuilletage $\xi F$ est d\'efini par
le champ de vecteurs~$Y$ dont le produit int\'erieur avec~$\omega$ vaut la
$1$-forme induite par $\alpha$ sur~$F$. Une singularit\'e de $\xi F$ est un
point de contact entre $\xi$ et~$F$ et a donc un signe : elle est positive ou
n\'egative selon que les orientations de $\xi$ et~$F$ co\"{\i}ncident ou non. Ce
signe est aussi celui de la divergence de~$Y$.

\medskip

Sur une vari\'et\'e de contact $(V,\xi)$, la projection $TV \to TV/\xi$ induit
un hom\'eomorphisme entre les champs de vecteurs sur~$V$ qui pr\'eservent~$\xi$
et les sections du fibr\'e quotient $TV/\xi \to V$. Autrement dit, si $\xi$ est
le noyau d'une $1$-forme~$\alpha$ et si $v$ est une fonction quelconque $V \to
\R$, il y a un unique champ de vecteurs $\nabla_\alpha v$ qui pr\'eserve~$\xi$
et v\'erifie $\alpha (\nabla_\alpha v) = v$. La flexibilit\'e des isotopies de
contact qui en r\'esulte fait que, comme en topologie diff\'erentielle, le
groupe $\DD (V,\partial V ;\xi)$ participe \`a des fibrations naturelles qui
peuvent aider \`a comprendre son type d'homotopie. Par exemple, si $F \subset V$
est une surface close et si $\PP(F,V ;\xi)$ d\'esigne l'espace des plongements de
$F$ dans~$V$ qui induisent sur~$F$ le m\^eme feuilletage caract\'eristique que
l'inclusion, la restriction \`a~$F$ d\'etermine une fibration
$$ \DD (V, \partial V ; \xi) \longrightarrow \PP (F, V ; \xi) \,. $$
D'autre part, en raison de la stabilit\'e des structures de contact, le groupe
$\DD (V,\partial V ;\xi)$ poss\`ede, dans $\DD (V,\partial V)$, un voisinage
qui se r\'etracte sur lui. La th\'eorie des surfaces convexes --~dont on rappelle
ci-dessous quelques points utiles~-- fournit de m\^eme, pour $\PP(F,V ;\xi)$, un
\'equivalent homotopique plus maniable (lemme~7).

\medskip

Soit $F$ une surface compacte, orient\'ee, proprement plong\'ee dans $(V,\xi)$
et \`a bord vide ou legendrien\footnote
{Une courbe legendrienne est une courbe int\'egrale de~$\xi$.}.
On dit que $F$ est {\it convexe\/} si elle a un voisinage produit
$$ U = F \times \R \supset F = F \times \{0\}, \quad \text{avec} \quad
   \partial U = \partial F \times \R \subset \partial V, $$
dans lequel l'action de~$\R$ par translations --~action engendr\'ee par le champ
de vecteurs $\partial_t$, $t \in \R$~-- pr\'eserve~$\xi$. Un tel voisinage~$U$
est dit {\it homog\`ene}.

Une surface convexe poss\`ede non seulement un voisinage homog\`ene mais tout un
syst\`eme fondamental. En effet, si $\nu$ est la section de $TV/\xi$ associ\'ee
au champ de vecteurs $\partial_t$ dans un voisinage homog\`ene $U = F \times
\R$, le fait que $\partial F$ soit legendrien assure que le champ de vecteurs
associ\'e \`a toute section du type $f \nu$, o\`u $f$ est une fonction qui ne
d\'epend que de~$t$, reste tangent \`a~$\partial V$. Cela r\'esulte, pour $\xi =
\ker \alpha$ et $v = \alpha(\nu)$, de la formule
$$ \nabla_\alpha (fv) - f \nabla_\alpha v = v Y_f \eqno(*) $$
o\`u $Y_f$ est le champ de vecteurs legendrien $d\alpha$-dual de $-df \res \xi$.

\smallskip

La convexit\'e est une propri\'et\'e tr\`es facile \`a d\'eceler. Une premi\`ere
condition n\'ecessaire est que chaque composante~$L$ de $\partial F$ poss\`ede,
sur $\partial V$, un ``voisinage homog\`ene'' ; on appelle ainsi, par extension,
un voisinage annulaire
$$ L \times \R \supset l = L \times \{0\} $$
dans lequel $\xi F$ est $\R$-invariant (figure~1). Dans ce voisinage, $\xi F$
est dirig\'e par un champ de vecteurs ind\'ependant de $t \in \R$ et plus
pr\'ecis\'ement du type $\cos(n\theta) \, \partial_\theta$, o\`u $\theta
\in \S^1$ param\`etre~$L$ et $n$~est un entier. Lorsque le bord de~$V$ est
lui-m\^eme convexe (ce qu'on supposera ensuite), un tel voisinage de~$L$ ne peut
exister que si $n$ n'est pas nul, {\it i.e.} si $\xi$ fait un nombre non nul de
tours autour de $\partial V$ le long de~$L$ dans le sens des aiguilles d'une
montre. On d\'ecrit maintenant les propri\'et\'es du feuilletage
caract\'eristique $\xi F$ li\'ees \`a la convexit\'e.

\begin{figure}[htbp]
\centering
\includegraphics[width=.5\textwidth]{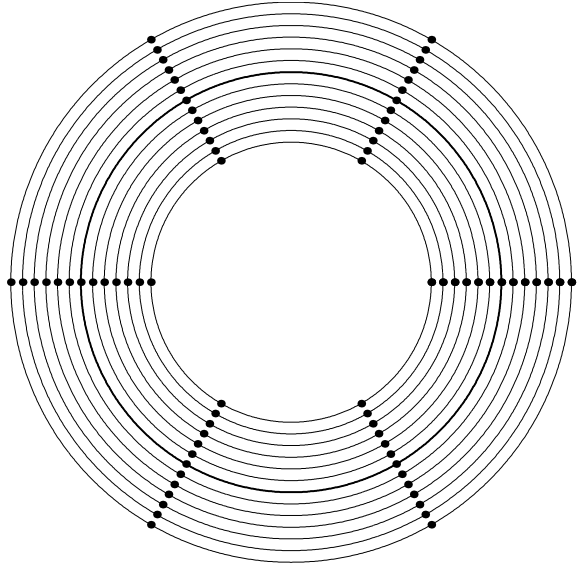}
\caption{Voisinage homog\`eme d'une courbe legendrienne~$L$ sur $\partial V$ :
dans cet exemple, l'entier~$n$ vaut~$3$.}
\end{figure}

\smallskip

On consid\`ere sur~$F$ un feuilletage singulier~$\sigma$ --~ensemble des courbes
int\'egrales orient\'ees d'un champ de vecteurs~-- et une multi-courbe~$\Gamma$ 
--~union finie de courbes ferm\'ees, simples et disjointes. On dit que $\Gamma$
{\it scinde\/}~$\sigma$ si $\sigma$ est transversal \`a~$\Gamma$ et est port\'e
par un champ de vecteurs qui, sur chaque composante de $F \setminus \Gamma$, ou
bien dilate l'aire et sort le long de~$\Gamma$, ou bien contracte l'aire et
rentre le long de~$\Gamma$. Une telle courbe, si elle existe et si $\sigma$ est
tangent \`a $\partial F$, est unique \`a isotopie pr\`es parmi les courbes qui
scindent.

\medskip

Les deux lemmes qui suivent sont d\'emontr\'es dans~\cite{Gi1} pour les surfaces
closes ; une fois les d\'efinitions pos\'ees, les d\'emonstrations s'\'etendent
sans difficult\'e au cas des surfaces \`a bord legendrien.

\begin{lemme} 
Soit $F$ une surface compacte, orient\'ee, proprement plong\'ee dans $(V,\xi)$
et \`a bord legendrien. Pour que $F$ soit convexe, il faut et il suffit que
$\partial F$ poss\`ede un voisinage homog\`ene sur $\partial V$ et que $\xi F$
soit scind\'e par une multi-courbe. De plus, cette seconde condition est remplie
d\`es que les propri\'et\'es suivantes sont satisfaites\up{ :}
\begin{itemize}
\item
chaque demi-orbite de $\xi F$ a pour ensemble limite une singularit\'e ou une
orbite ferm\'ee\up{ ;}
\item
les orbites ferm\'ees de $\xi F$ sont toutes non d\'eg\'en\'er\'ees --~i.e.
hyperboliques\up{ ;}
\item
aucune orbite de $\xi F$ ne va d'une singularit\'e n\'egative \`a une
singularit\'e positive.
\end{itemize}
\end{lemme}

Comme on l'explique dans~\cite{Gi1}, cette caract\'erisation montre que, si le
bord de~$F$ est legendrien et poss\`ede un voisinage homog\`ene sur $\partial V$,
on peut rendre~$F$ convexe par une isotopie relative au bord et arbitrairement
petite. Par ailleurs, quand $F$ est convexe, le choix d'un voisinage homog\`ene
$U = F \times \R$ d\'etermine une multi-courbe $\Gamma_{\!U}$ qui scinde $\xi
F$ : c'est l'ensemble des points de~$F$ o\`u le vecteur $\partial_t$, $t \in \R$,
appartient au plan~$\xi$. Parmi toutes les multi-courbes qui scindent $\xi F$,
celles qui sont ainsi associ\'ees aux voisinages homog\`enes de~$F$ ont la
particularit\'e d'aboutir aux points de contact entre $\xi$ et~$\partial V$.

\begin{lemme} 
Soit $F \subset (V,\xi)$ une surface convexe, soit $U$ un voisinage homog\`ene
de~$F$ et soit $\Gamma_{\!U} \subset F$ la multi-courbe associ\'ee \`a~$U$.

\alinea{\textbf{a)}}
L'espace $\FF (F,\Gamma)$ des feuilletages singuliers de~$F$ qui sont
scind\'es par~$\Gamma_{\!U}$ et tangents \`a $\partial F$ est contractile.

\alinea{\textbf{b)}}
Soit $\PP(F,V)$ l'espace des plongements de~$F$ dans~$V$ qui co\"\i ncident avec
l'inclusion sur $\partial F$. Il existe une application continue
\begin{align*}
   \FF (F, \Gamma_{\!U}) & \longrightarrow \PP (F,V) \\
   \sigma & \longmapsto \psi_\sigma
\end{align*}
ayant les propri\'et\'es suivantes\up{ :}
\begin{itemize}
\item
$\psi_{\xi F}$ est l'inclusion\up{ ;}
\item
pour tout $\sigma$, la surface $\psi_\sigma (F)$ est incluse dans~$U = F \times
\R$ et transversale au champ de vecteurs $\partial_t$, $t \in \R$\up{ ;}
\item
pour tout~$\sigma$, le feuilletage caract\'eristique de $\psi_\sigma (F)$ n'est
autre que $\psi_{\sigma*} \sigma$.
\end{itemize}
\end{lemme}

Sous les hypoth\`eses du lemme ci-dessus et avec ses notations, on dira qu'un
plongement quelconque $\psi \in \PP(F,V)$ est {\it adapt\'e \`a $\Gamma_{\!U}$}
si le feuilletage caract\'eristique de $\psi(F)$ est scind\'e par $\psi
(\Gamma_{\!U})$.

\begin{lemme} 
Soit $F \subset (V,\xi)$ une surface convexe, $\Gamma$ la multi-courbe sur~$F$
associ\'ee \`a un voisinage homog\`ene quelconque, $\PP(F,V;\Gamma)$ l'espace
des plongements de~$F$ dans~$V$ adapt\'es \`a~$\Gamma$ et $\PP(F,V;\xi)$ le
sous-espace de ceux qui induisent le m\^eme feuilletage caract\'eristique que
l'inclusion. Le plongement canonique de $\PP(F,V;\xi)$ dans $\PP(F,V;\Gamma)$
est une \'equivalence d'homotopie.
\end{lemme}

\begin{proof}
Vu le lemme~6, il suffit d'associer contin\^ument, \`a tout plongement $\psi \in
\PP(F,V ;\Gamma)$, un champ de vecteurs de contact $X(\psi)$ sur $(V,\xi)$ dont
le flot soit transversal \`a $\psi(F)$ et plonge $\psi(F) \times \R$ dans~$V$ en
envoyant $\partial \psi(F) \times \R$ dans $\partial V$. On fixe une m\'etrique
riemannienne et on choisit une fonction continue positive $d(\psi)$ telle que,
pour tout $\psi \in \PP(F,V)$, le $d(\psi)$-voisinage de~$\psi(F)$ soit un tube.
On note $\XX (\psi)$, $\psi \in \PP(F,V ;\Gamma)$, l'espace des champs de
vecteurs de contact transversaux \`a $\psi(F)$, tangents \`a $\partial V$ et
dont le support est confin\'e \`a l'int\'erieur du $d(\psi)$-voisinage de~$\psi
(F)$. Ainsi, le flot de tout champ $X \in \XX (\psi)$ plonge $\psi(F)
\times \R$ dans~$V$ et envoie $\partial \psi(F) \times \R$ dans $\partial V$. On
peut alors pr\'elever contin\^ument un champ $X(\psi)$ dans $\XX (\psi)$
pour les raisons suivantes :
\begin{itemize}
\item
chaque espace $\XX (\psi)$ est contractile --~et non vide ;
\item
tout champ $X \in \XX (\psi)$ appartient aussi \`a $\XX (\psi')$
pour $\psi'$ assez proche de~$\psi$.
\end{itemize}
\end{proof}

Les consid\'erations ci-devant sugg\`erent qu'il est plus commode d'\'etudier la
topologie de l'espace $\PP(F,V ;\xi)$ --~et ainsi du groupe $\DD(V,\partial V ;
\xi)$~-- lorsque $F$ est une surface convexe. Par ailleurs, cette propri\'et\'e
est g\'en\'eriquement satisfaite pour peu que $F$ soit close ou que $\partial F$
soit une courbe legendrienne et poss\`ede un voisinage homog\`ene sur $\partial
V$. Cette condition n'est bien s\^ur pas toujours remplie : il se peut m\^eme que
$\partial F$ ne soit isotope \`a aucune courbe legendrienne sur $\partial V$.
Cependant la convexit\'e --~cette fois de $\partial V$~-- permet de contourner
ce probl\`eme.

\smallskip

On suppose d\'esormais que $(V,\xi)$ est une vari\'et\'e de contact dont le bord
est convexe ; comme $\partial V$ n'a pas de bord, cette hypoth\`ese n'est gu\`ere
restrictive. On choisit un voisinage collier homog\`ene
$$ W = \partial V \times [0,\infty \cl[ \supset \partial V = \partial V \times \{0\} $$
et on note $\Delta$ la multi-courbe associ\'ee sur $\partial V$. On consid\`ere
alors l'espace $\SC(V,\Delta)$ des structures de contact sur~$V$ pour lesquelles
$\partial V$ est convexe et a un feuilletage caract\'eristique scind\'e par~$
\Delta$. Le point est que la composante connexe de~$\xi$ dans $\SC (V,\Delta)$
contient une structure de contact pour laquelle $\partial F$ est isotope \`a une
courbe legendrienne ayant un voisinage homog\`ene sur~$\partial V$. C'est ce que
garantit le lemme~6 : en effet, si $C \subset \partial V$ est une courbe isotope
\`a $\partial F$ dont chaque composante rencontre~$\Delta$, il existe sur
$\partial V$ un feuilletage singulier scind\'e par~$\Delta$ et qui, sur un
voisinage $C \times \R$ de $C = C \times \{0\}$, est tangent aux courbes $C
\times \{t\}$, $t \in \R$ \ (voir~[Gi1, exemple II.3.7] pour plus de d\'etails).
Ce fait est compl\'et\'e par le r\'esultat suivant :

\begin{proposition} 
Le type d'homotopie du groupe $\DD (V,\partial V ; \xi)$ ne d\'epend que de la
composante connexe de $\SC(V,\Delta)$ qui contient~$\xi$.
\end{proposition}

La d\'emonstration utilise les notations et le lemme ci-dessous. Pour toute
fonction $d \from \partial V \to \op] 0, \infty \cl[$, on pose
$$ W^d = \bigl\{ (p,t) \in W = \partial V \times [0,\infty \cl[ \mid
   t \le d(p) \bigr\} \,. $$
On d\'esigne de plus par $\PP(W^d,V ;\xi)$ l'espace des plongements de contact de
$(W^d,\xi)$ dans $(V,\xi)$ dont la restriction \`a $\partial V$ est l'inclusion.

\begin{lemme} 
Quelle que soit la fonction $d \from \partial V \to [0,\infty \cl[$, l'espace
$\PP(W^d,V ;\xi)$ est contractile.
\end{lemme}

\begin{proof}
Pour tout plongement $\psi \in \PP(W^d,V ;\xi)$, le champ de vecteurs $\partial_t
\psi$ --~o\`u $t$ est la coordonn\'ee dans $\op[0,\infty\cl[$~-- se prolonge
\`a~$V$ en un champ de vecteurs $X(\psi)$ qui pr\'eserve~$\xi$ et varie
contin\^ument avec~$\psi$. On consid\`ere alors les champs de vecteurs de
contact
$$ X_s = (1-s) X(\psi_0) + s X(\psi), \qquad s \in [0,1], $$
o\`u $\psi_0$ d\'esigne l'inclusion $W^d \to V$. L'int\'egration du champ $X_s$
\`a partir de chaque point $p \in \partial V$ pendant un temps~$d(p)$ fournit un
plongement $\psi_s$ de $W^d$ dans~$V$. Mieux, $\psi_s$ est un plongement de
contact de $(W^d,\xi)$ dans $(V,\xi)$. En effet, comme $\psi$ pr\'eserve~$\xi$
et induit l'identit\'e sur $\partial V$, les sections de $TV/\xi$ d\'efinies par
$X(\psi)$ et $X(\psi_0)$ co\"{\i}ncident en tout point de $\partial V$ o\`u le
plan~$\xi$ est transversal \`a $\partial V$, donc sur un ouvert dense de
$\partial V$. Il en r\'esulte que les champs $X_s$ induisent tous la m\^eme
section de $TV/\xi$ le long de $\partial V$, ce qui implique la propri\'et\'e
voulue. Les plongements $\psi_s$ forment ainsi un chemin canonique joignant
$\psi_0$ \`a~$\psi$ dans $\PP (W^d,V ;\xi)T$.
\end{proof}

\begin{proof}[D\'emonstration de la proposition~8]
Pour toute fonction positive~$d$ sur $\partial V$, la restriction \`a~$W^d$
d\'efinit une application
$$ \DD (V,\partial V ;\xi) \longrightarrow \PP (W^d,V ;\xi) $$
qui, vu les crit\`eres g\'en\'eraux de J.~Cerf~\cite{Ce} et la flexibilit\'e des
isotopies de contact, est une fibration localement triviale. D'apr\`es le lemme
9, le groupe $\DD (V,\partial V ;\xi)$ se r\'etracte donc sur son sous-groupe
$\DD (V,W^d ;\xi)$ form\'e des diff\'eomorphismes qui valent l'identit\'e sur
$W^d$.

Soit maintenant $\xi_s$, $s \in [0,1]$, un chemin dans $\SC (V,\Delta)$ partant
du point $\xi_0 = \xi$. Le lemme~6 et la m\'ethode de J.~Moser fournissent un
plongement de contact $\phi$ de $(V,\xi_1)$ dans $(V,\xi)$ qui envoie $\partial
V$ dans~$W$ sur une surface transversale au champ de vecteurs $\partial_t$.
Autrement dit, $\phi (\partial V)$ est, dans $W = F \times [0,\infty\cl[$, le
graphe d'une fonction $d \from \partial V \to \op]0,\infty\cl[$. D'autre part,
l'ouvert $W_1 = \phi^{-1}(W)$ est un voisinage collier homog\`ene de $\partial
V$ dans $(V,\xi_1)$. De plus, pour toute fonction $e \from \partial V \to \op]0,
\infty\cl[$,
$$ W_1^e  := (W_1)^e = \phi^{-1} (W^{d+e}) \,. $$
Par cons\'equent, les groupes $\DD(V,\partial V;\xi)$ et $\DD(V,\partial V;
\xi_1)$ se r\'etractent respectivement sur les sous-groupes $\DD(V,W^{d+e};\xi)$
et $\DD(V,W_1^e;\xi_1)$, sous-groupes entre lesquels $\phi$ induit clairement un
hom\'eomorphisme.
\end{proof}

\section{Voisinages r\'etractiles et isotopies de contact}

Dans cette partie, $(V,\xi)$ d\'esigne toujours une vari\'et\'e de contact de
dimension~$3$, compacte et \`a bord convexe. On se donne en outre dans $(V,\xi)$
une surface $F$, proprement plong\'ee et convexe, et on note $\Gamma \subset F$
la multi-courbe associ\'ee \`a un voisinage homog\`ene quelconque. Compte tenu
des remarques g\'en\'erales de la partie~A, on \'etudie maintenant l'espace $\PP
(F,V ;\Gamma)$ des plongements adapt\'es \`a~$\Gamma$. Cette \'etude s'inspire du
travail de V.~Colin~\cite{Co} sur les plongements des sph\`eres. La cl\'e de la
d\'emonstration des th\'eor\`emes principaux est ainsi le r\'esultat suivant :

\begin{proposition} 
Soit $\psi \in \PP (F,V ;\Gamma)$, soit $\CC$ la composante connexe de
$\psi$ dans $\PP (F,V)$ et soit $\CC(\Gamma)$ l'intersection de $\CC$ avec
$\PP(F,V ;\Gamma)$. Pour que $\CC(\Gamma)$ soit connexe, il suffit
qu'il connecte~$\psi$ \`a tout plongement $\psi' \in \CC (\Gamma)$ qui a
les propri\'et\'es suivantes\up{ :}
\begin{itemize}
\item
$\psi'(F) \cap \psi(F) = \partial F$\up{ ;}
\item
$\psi'(F) \cup \psi(F)$ borde dans~$V$ un domaine hom\'eomorphe au produit $F
\times \D^1$.
\end{itemize}
\end{proposition}

Gr\^ace \`a cette proposition, la connaissance des structures de contact sur le
produit $F \times [0,1]$ permet, dans certains cas, de prouver que $\CC
(\Gamma)$ est connexe. \`A titre d'exemple, on reproduit l'argument donn\'e dans
\cite{Co} pour traiter le cas des sph\`eres. On rappelle qu'une vari\'et\'e de
contact $(V,\xi)$ est dite {\it tendue\/} si aucun disque plong\'e dans~$V$
n'est tangent \`a~$\xi$ en tous les points de son bord.

\begin{corollaire}[V.~Colin] 
Si la vari\'et\'e de contact $(V,\xi)$ est tendue et si la surface~$F$ est une
sph\`ere, l'espace $\CC (\Gamma)$ est connexe.
\end{corollaire}

\begin{proof}
Soit $\eta$ une structure de contact sur $W = \S^2 \times [0,1]$. Si $(W,\eta)$
est tendue et \`a bord convexe, un th\'eor\`eme de Y.~Eliashberg~\cite{El} (voir
aussi [Gi2, lemme 2.17]) assure que $\eta$ est isotope, relativement au
bord, \`a une structure de contact pour laquelle toutes les sph\`eres $\S^2
\times \{*\}$ sont convexes. Joint \`a la proposition~10, ce fait entra\^\i ne
imm\'ediatement le corollaire.
\end{proof}

La d\'emonstration de la proposition~10 repose sur la notion de voisinages
r\'etractiles.

\begin{definition} 
Soit $F_0$ une surface compacte proprement plong\'ee dans $(V,\xi)$. On dira
qu'un voisinage $U_0$ de~$F_0$ est {\it r\'etractile\/} si, pour tout voisinage
$U$ de~$F_0$, il existe une isotopie de contact $\phi_s \from V \to V$, $s \in
[0,1]$, ayant les propri\'et\'es suivantes :
\begin{itemize}
\item
$\phi_0 = \id$ ;
\item
$\phi_s \res {F_0} = \id$ pour tout $s \in [0,1]$ ;
\item
$\phi_1 (U_0) \subset U$.
\end{itemize}
\end{definition}

L'exemple typique de voisinages r\'etractiles est le suivant :

\begin{lemme} 
Soit $U_0 = F_0 \times [-1,1]$ un voisinage de $F_0 = F_0 \times \{0\}$. Si les
surfaces $F_t = F \times \{t\}$, $t \ne 0$, sont toutes convexes, $U_0$ est un
voisinage r\'etractile de~$F_0$.
\end{lemme}

\begin{proof}
On fabrique, pour $\eps>0$ donn\'e, une isotopie de contact qui r\'etracte $U_0$
dans $F_0 \times [-\eps,\eps]$. Pour tout $t \ne 0$, comme $F_t$ est convexe, il
existe un champ de vecteurs de contact~$X_t$ transversal \`a~$F_t$ et tangent
\`a $\partial V$. Quitte \`a changer~$X_t$ en son oppos\'e, on suppose que $X_t$
pointe vers~$F_0$ le long de~$F_t$. On se donne en outre un voisinage $J_t$ de
$t$ dans $[-1,1] \setminus \{0\}$ tel que chaque champ~$X_t$ reste transversal
\`a $F_{t'}$ pour tout $t' \in J_t$. On extrait ensuite de la famille $J_t$ un
recouvrement fini $J_i = J_{t_i}$ de $[-1,-\eps] \cup [\eps,1]$ et on prend une
partition de l'unit\'e $f_i$ subordonn\'ee \`a ce recouvrement. Le champ de
vecteurs $\Ol X = \sum f_i X_{t_i}$ ne pr\'eserve pas $\xi$ mais, d'apr\`es la
formule~$(*)$, il diff\`ere d'un champ de vecteurs de contact par un champ
legendrien~$Y$ nul pr\`es de~$F_0$ et tangent aux surfaces $F_t$ ainsi qu'\`a
$\partial V$. Il suffit alors d'int\'egrer le champ $\Ol X + Y$ pendant un temps
assez long pour obtenir l'isotopie de contact voulue.
\end{proof}

On d\'emontre maintenant la proposition~10. L'ingr\'edient essentiel est le
lemme ci-dessous qu'on admet provisoirement :

\begin{lemme} 
Tout chemin dans $\PP (F,V)$ qui relie deux points de $\PP (F,V ;\Gamma)$ est
homotope, relativement \`a ses extr\'emit\'es, \`a un chemin de plongements dont
les images ont chacune un voisinage r\'etractile.
\end{lemme}

\begin{proof}[D\'emonstration de la proposition~10]
Soit $\psi_s$, $s \in [0,1]$, un chemin dans $\CC$ qui relie un point
$\psi_0 \in \CC (\Gamma)$ \`a $\psi_1 = \psi$. On construit ci-dessous un
chemin $\psi'_s$ dans $\CC (\Gamma)$ qui joint $\psi'_0 = \psi_0$ \`a un
plongement $\psi'_1$ dont l'image est disjointe de celle de~$\psi_1$ et borde
avec celle-ci un produit. 

Pour tout $s \in [0,1]$, on pose $F_s = \psi_s(F)$. Gr\^ace au lemme~14, il est
possible de trouver des points $s_0 = 0 < s_1 < \dots < s_k = 1$ tels que, pour
$0 \le i \le k-1$, les surfaces $F_s$, $s \in [s_i,s_{i+1}]$ soient toutes
incluses dans un
voisinage r\'etractile $U_i$ de $F_{s_i}$. Comme $F_0$ est convexe, on prend en
fait pour $U_0$ un voisinage produit $F_0 \times [-1,1]$ dans lequel toutes les
surfaces $F_0 \times \{t\}$ sont convexes. Si $F$ est close, on pose
$$ F'_s = F_0 \times \Bigl\{ \tfrac s {s_1} \Bigr\} \subset U_0 \quad
   \text{pour $s \in [0,s_1]$.} $$
Si $F$ a un bord, on se donne une fonction $f \from F \to [0,1]$ qui vaut $0$
sur $\partial F$ et $1$ hors du $\eps$-voisinage de $\partial F$ ; pour $\eps$
assez petit, les graphes des fonctions $r f$, $r \in [-1,1]$, sont convexes et
couvrent la r\'eunion des surfaces $F_s$, $s \in [0,s_1]$. On pose alors
$$ F'_s = \Bigl\{ \bigl( p, \tfrac s {s_1} f(p) \bigr), \ p \in F \Bigr\}
   \subset U_0 \quad
   \text{pour $s \in [0,s_1]$.} $$
Comme les surfaces $F'_s$ sont convexes, il est possible de les param\'etrer par
des plongements $\psi'_s \in \CC (\Gamma)$.

Maintenant, par hypoth\`ese, $F_{s_1}$ est contenue dans $U_0$ et $U_1$ est un
voisinage r\'etractile de $F_{s_1}$. Il existe donc une isotopie de contact
$\phi_r^1$, $r \in [0,1]$, qui r\'etracte $U_1$ dans~$U_0$. Quitte \`a changer
le param\'etrage en~$r$, on suppose que $\phi_r^1$, pour tout $r \in [0,1]$,
envoie $F_{(1-r)s_1+rs_2}$ dans~$U_0$. On pose alors
\begin{gather*}
   F'_s = \bigl( \phi_r^1 \bigr)^{-1} \bigl( F'_{s_1} \bigr)
   \quad \text{pour $s \in [s_1,s_2]$,} \\
   \text{o\`u} \quad r = \frac{ s-s_1 }{ s_2-s_1 } .
\end{gather*}
Par construction, chaque surface $F'_s$ est disjointe de~$F_s$ et borde avec
elle un produit. De plus, $F'_s$ est convexe et se laisse donc param\'etrer par
un plongement $\psi'_s \in \CC (\Gamma)$. On continue en observant que
$F_{s_2}$ est incluse dans $U^1 = (\phi_1^1)^{-1} (U_0)$. On r\'etracte $U_2$
dans~$U^1$ par une isotopie de contact $\phi_r^2$ qui, pour tout $r \in [0,1]$,
envoie $F_{(1-r)s_2+rs_3}$ dans $U^1$. On d\'efinit alors $F'_s$ par une formule
du m\^eme type que ci-dessus et on obtient le chemin $\psi'_s$ tout entier en
r\'ep\'etant $k$~fois l'op\'eration.
\end{proof}

\begin{proof}[D\'emonstration du lemme~14]
Soit $\psi_s$, $s \in [0,1]$, un chemin dans $\PP(F,V)$. On peut trouver des
points $s_0 = 0 < s_1 < \dots < s_k = 1$ tels que, pour $0 \le i \le k-1$, les
surfaces $F_s$, $s \in [s_i,s_{i+1}]$, soient toutes contenues dans un domaine
produit\footnote
{\`A strictement parler, il faut consid\'erer ici des {\it produits pinc\'es\/}
$U_i \simeq F \times_\partial [-1,1]$, o\`u le terme~\cite{Ha2} et la notation
d\'esignent le quotient de $F \times [-1,1]$ par la relation qui \'ecrase sur un
point chaque segment $\{p\} \times [-1,1]$, $p \in \partial F$. On laisse au
lecteur le soin de faire les quelques ajustements n\'ecessaires.}
$U_i \simeq F \times [-1,1]$ plong\'e dans~$V$. On note alors~$F^i$ une des deux
copies de~$F$ qui forment $\partial U_i$. Maintenant, pour $s \in[s_i,s_{i+1}]$,
on remplace l'isotopie $\psi_s$ par une isotopie $\psi'_s$ dont les images, pour
$s \le (s_i+s_{i+1})/2$ (resp. $s \ge (s_i+s_{i+1})/2$), balaient le produit que
bordent $\psi_{s_i}(F)$ et $F^i$ (resp. $F^i$ et $\psi_{s_{i+1}}(F)$). Gr\^ace
\`a cette astuce emprunt\'ee \`a V.~Colin~\cite{Co} et au lemme~13, il suffit
d'\'etablir le r\'esultat suivant :
\def\qed{}
\end{proof}

\begin{lemme} 
Soit $\xi$ une structure de contact sur $W = F \times [-1,1]$. Si le bord de
$(W,\xi)$ est convexe, il existe une isotopie relative au bord $\phi_s \from W
\to W$, $s \in [0,1]$, telle que toutes les surfaces $\phi_1 (F \times \{t\})$
soient convexes, sauf un nombre fini qui portent un feuilletage caract\'eristique
ayant les propri\'et\'es suivantes\up{ :}
\begin{itemize}
\item
les singularit\'es et les orbites ferm\'ees sont toutes non d\'eg\'en\'er\'ees
\up(c'est-\`a-dire hyperboliques\up{) ;}
\item
l'ensemble limite de chaque demi-orbite est une singularit\'e ou une orbite
ferm\'ee\up{ ;}
\item
une orbite et une seule va d'une selle n\'egative vers une selle positive.
\end{itemize}
\end{lemme}

La d\'emonstration de ce lemme repose sur les techniques de~\cite{Gi2} et passe par
une construction explicite d'isotopies :

\begin{lemme} 
On consid\`ere sur $\R^3$ la structure de contact d'\'equation $dz + x\,dy = 0$
et on se donne trois nombres $a, b, c \in \op]0, 1\cl[$. Il existe une fonction
$f \from \R^3 \to \R$, \`a support dans $[1-a,1+a] \times [b,b] \times [-c,c]$,
ayant les propri\'et\'es suivantes\up{ :}
\begin{itemize}
\item
$f(x,0,z) = 0$ pour tout $x, z \in \R$\up{ ;}
\item
les transformations $\phi_s \from \R^3 \to \R^3$, $s \in [0,1]$, d\'efinies par
$$ \phi_s (x,y,z) = \bigl( x, y, z + s f(x,y,z) \bigr) $$
forment une isotopie\up{ ;}
\item
pour tout $s$ et tout $z$, le feuilletage caract\'eristique de $\phi_s (\R^2
\times \{z\})$ n'a aucune singularit\'e hors de l'axe $\{ y=0 \}$, en a au plus
deux sur cet axe et exactement deux si $s=1$ et $z=0$.
\end{itemize}
\end{lemme}

\begin{proof}
On pose $f(x,y,z) = - u(x) v(y) w(z)$ o\`u $u$, $v$ et~$w$ sont des fonctions
qui remplissent les conditions ci-dessous et sont nulles pr\`es du bord de leurs
domaines respectifs :
\begin{itemize}
\item
$u \from [1-a,1+a] \to [0,2]$ a un graphe qui, vu dans $\R^2$, rencontre la
diagonale en deux points et est strictement croissant-et-concave entre eux ;
\item
$u \from [-b,b] \to [-c/4,c/4]$ s'annule en $0$ avec d\'eriv\'ee $1$ mais est
non nulle et a une d\'eriv\'ee strictement inf\'erieure \`a~$1$ en tout autre
point \`a l'int\'erieur de son support ;
\item
$w \from [-c,c] \to [0,1]$ vaut~$1$ pr\`es de~$0$ et sa d\'eriv\'ee est partout
born\'ee par $2/c$ en valeur absolue.
\end{itemize}
Le feuilletage caract\'eristique de la surface $\phi_s (\R^2 \times \{z\})$ est
alors engendr\'e, \`a la source, par le champ de vecteurs
$$ \bigl( x - s u(x) v'(y) w(z) \bigr) \, \partial_x + s u'(x) v(y) w(z) \,
   \partial_y $$
et on v\'erifie sans peine que les propri\'et\'es souhait\'ees sont satisfaites.
\end{proof}

\begin{proof}[D\'emonstration du lemme~15]
Le lemme 2.10 et la remarque~2.11 de \cite{Gi2} permettent de supposer que les
valeurs de~$t$ pour lesquelles $F_t$ n'est pas convexe constituent un compact
d\'enombrable $\Sigma \subset [-1,1]$ et que, pour tout $t \in \Sigma$, le
feuilletage caract\'eristique $\xi F_t$ a les propri\'et\'es suivantes :
\begin{enumerate}
\item[1)]
les singularit\'es de $\xi F_t$ sont toutes non d\'eg\'en\'er\'ees ;
\item[2)]
chaque demi-orbite de $\xi F_t$ a pour ensemble limite une singularit\'e ou une
orbite ferm\'ee ;
\item[3)]
$\xi F_t$ a exactement une orbite ``exceptionnelle'' qui ou bien va d'une selle
n\'egative \`a une selle positive ou bien se referme avec une holonomie tangente
\`a l'identit\'e mais dont la d\'eriv\'ee seconde est non nulle ;
\item[4)]
$F_t$ est partag\'ee en deux sous-surfaces planaires $F_t^+$ et $F_t^-$ (non
n\'ecessairement connexes) et, le long du bord commun, $\xi F_t$ sort de $F_t^+$
pour entrer dans $F_t^-$.
\end{enumerate}
Dans ces conditions, les points $t_0 \in \Sigma$ pour lesquels $\xi F_{t_0}$ a
une orbite ferm\'ee d\'eg\'en\'er\'ee forment un ensemble $\Sigma_1$ fini. En
effet, vu la propri\'et\'e~4), cette orbite ferm\'ee~$C_{t_0}$ est enti\`erement
contenue dans $F_{t_0}^+$ ou dans $F_{t_0}^-$ et s\'epare ladite r\'egion. Par
suite, aucune orbite de $\xi F_{t_0}$ partant de $C_{t_0}$ ne revient vers
$C_{t_0}$. Compte tenu des propri\'et\'es 1) et 2), si un feuilletage $\xi F_t$
pr\'esente alors une orbite ferm\'ee d\'eg\'en\'er\'ee pour $t$ voisin de~$t_0$,
celle-ci reste dans un voisinage de $C_{t_0}$. Or le lemme 2.12 de~\cite{Gi2}
exclut cette \'eventualit\'e pour $t \ne t_0$.

Soit maintenant $t_0$ un point de $\Sigma_1$ et soit $A_{t_0}$ un voisinage de
$C_{t_0}$ dans lequel $\xi F_{t_0}$ est non singulier, transversal au bord et
sans autre orbite ferm\'ee que $C_{t_0}$. On se donne un nombre $\eps>0$ tel que
tout feuilletage $\xi F_t$, $t \in [t_0-\eps, t_0+\eps]$, satisfasse aux
conditions suivantes :
\begin{enumerate}
\item[5)]
$\xi F_t$ est transversal \`a $\partial A_{t_0}$ et \`a $\partial F_{t_0}^\pm$ \ 
(on identifie ici $F_t$ \`a $F_{t_0}$ via la projection sur~$F$) ;
\item[6)]
les singularit\'es de $\xi F_t$ sont toutes non d\'eg\'en\'er\'ees de m\^eme que
les orbites ferm\'ees situ\'ees hors de $A_{t_0}$.
\end{enumerate}
Soit $p_-$ et $p_+$ deux points de $C_{t_0}$. Sur un voisinage $U_\pm$ de
$p_\pm$, on peut trouver des coordonn\'ees $(x,y,z) \in \pm [1-a,1+a] \times
[-b,b] \times [-c,c]$ \ (o\`u $a$, $b$ et $c$ sont des r\'eels positifs petits)
dans lesquelles :
\begin{enumerate}
\item[7)]
la $1$-forme $\pm (dz - x \, dy)$ d\'efinit et cooriente $\xi$ sur $U_\pm$ ;
\item[8)]
le param\`etre $t$ est une fonction croissante de la seule coordonn\'ee $z$ (et
du signe $\pm$) qui vaut $t_0$ pour $z=0$ ;
\item[9)]
l'holonomie du feuilletage caract\'eristique $\xi F_{t_0}$ entre les segments
$T_- = \{-1+A\} \times [-b,b] \times \{0\} \subset U_-$ et $T_+ = \{1-a\} \times
[-b,b] \times \{0\} \subset U_+$ respecte la coordonn\'ee~$y$.
\end{enumerate}
On prend alors pour $\phi_s$, $s \in [0,1]$, une isotopie \`a support dans $U_+
\cup U_-$ et de la forme
$$ \phi_s (x,y,z) = \bigl( x, y, z \pm s f (\pm x, y, \pm z) \bigr), \qquad
   (x,y,z) \in U_\pm, $$
o\`u $f$ est une fonction fournie par le lemme~16. Par construction, le
feuilletage caract\'eristique $\xi \, \phi_1(F_{t_0})$ ne poss\`ede plus
d'orbite ferm\'ee dans $\phi_1(A_{t_0})$ mais pr\'esente une orbite qui va d'une
selle n\'egative \`a une selle positive. La propri\'et\'e~9) et la sym\'etrie de
l'isotopie $\phi_s$ permettent en outre de contr\^oler enti\`erement le type
topologique de ce feuilletage. Le lemme~2.12 de~\cite{Gi2} d\'ecrit alors les
feuilletages $\xi \, \phi_1(F_t)$ pour $t$ voisin de $t_0$ et assure notamment
que la modification effectu\'ee n'introduit aucune orbite ferm\'ee
d\'eg\'en\'er\'ee pour $t \in [t_0-\eps,t_0+\eps]$.
\end{proof}

\begin{figure}[htbp]
\centering
\includegraphics[width=.86\textwidth]{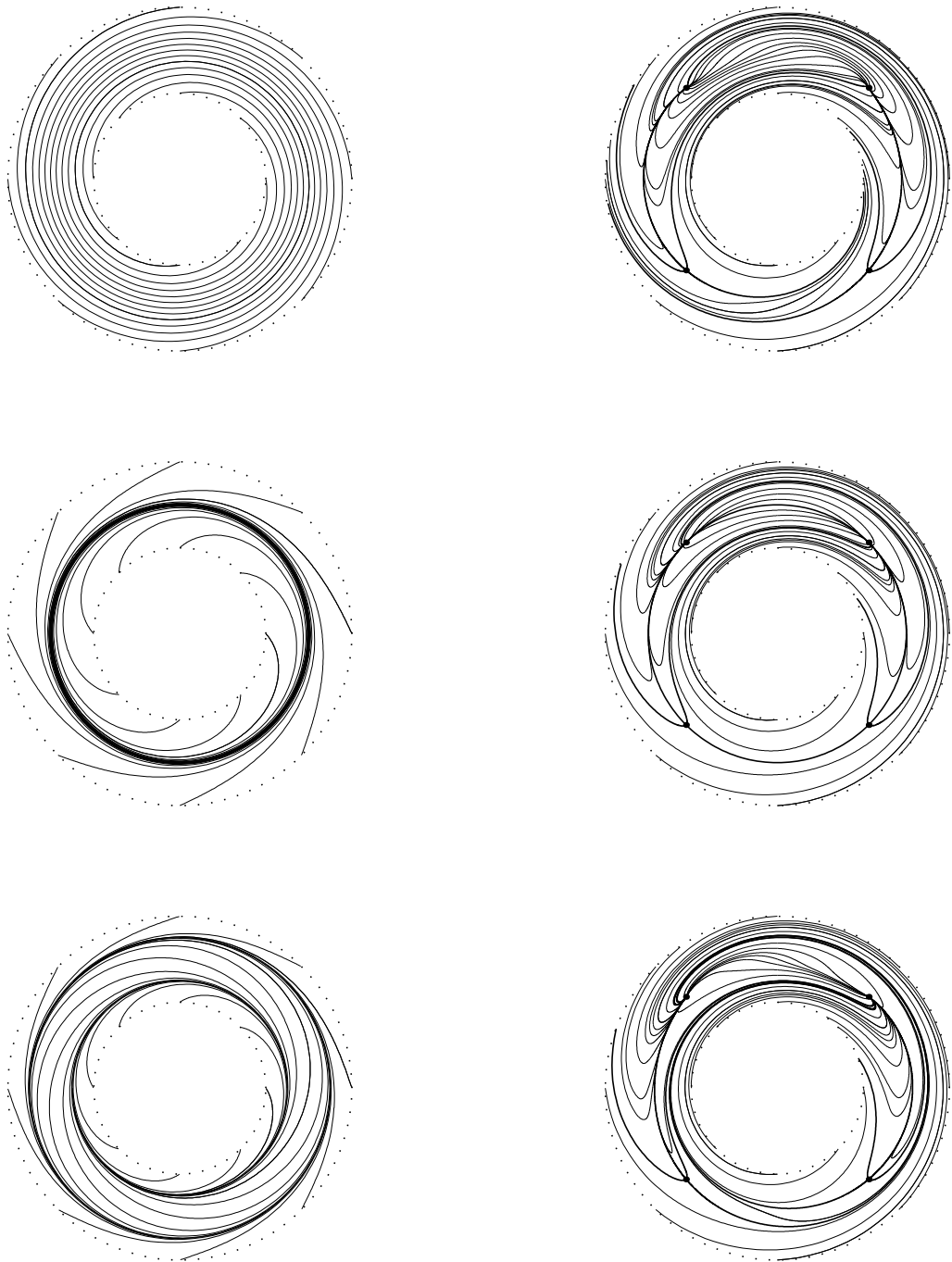}
\caption{
La partie gauche montre le feuilletage caract\'eristique de~$F_t$ dans l'anneau
$A_{t_0}$ pour $t<t_0$ (en haut), $t=t_0$ (au centre) et $t>t_0$ (en bas). La
partie droite affiche en regard le feuilletage de $\phi_1(F_t)$. Une bifurcation
organis\'ee autour d'une connexion de selles remplace la naissance de deux
orbites ferm\'ees.}
\end{figure}

\section{D\'emonstration des th\'eor\`emes}

Les th\'eor\`emes 1, 3 et~4 ont pour souche commune le r\'esultat suivant, d\^u
\`a Y.~Eliashberg~\cite{El} :

\begin{theoreme}[Y.~Eliashberg] 
Pour toute structure de contact tendue~$\xi$ sur la boule~$\D^3$, le groupe
$\DD (\D^3, \partial\D^3 ; \xi)$ est connexe --~et m\^eme contractile.
\end{theoreme}

\begin{proof}
D'apr\`es le th\'eor\`eme~2.4.2 de \cite{El}, l'espace $\SCT (\D^3 ;\xi)$ des
structures de contact tendues sur~$\D^3$ qui co\"{\i}ncident avec~$\xi$ le long
du bord est contractile. D'autre part, le groupe $\DD (\D^3,\partial \D^3)$
agit sur $\SCT (\D^3 ;\xi)$ et la m\'ethode de J.~Moser montre que cette action, 
restreinte au point~$\xi$, d\'efinit une fibration
$$ \DD (\D^3,\partial \D^3) \longrightarrow \SCT (\D^3 ;\xi) $$
dont la fibre n'est autre que $\DD (\D^3,\partial \D^3 ;\xi)$. Le fait que ce
groupe soit connexe (resp. contractile) d\'ecoule alors du th\'eor\`eme de
J.~Cerf~\cite{Ce} (resp. du th\'eor\`eme de A.~Hatcher~\cite{Ha1}) qui affirme
que $\DD (\D^3,\partial \D^3)$ est connexe (resp. contractile).
\end{proof}

On d\'esigne d\'esormais par~$S$ une surface compacte, connexe et orien\-table,
par $\pi \from V \to S$ le fibr\'e des \'el\'ements de contact coorient\'es sur
$S$ et par $\xi = \xi_S$ la structure de contact canonique sur~$V$.

\subsection*
{D\'emonstration du th\'eor\`eme~1 : cas du groupe $\DD (V,\partial V ;\xi)$}

On suppose d'abord que $S$ est le disque unit\'e $\D^2$ de $\R^2$ et on assimile
$V$ au tore plein $\D^2 \times \R/2\pi\Z$ en rep\'erant chaque droite tangente
\`a~$S$ en un point $(x_1,x_2)$ par l'angle~$\theta$ de sa normale directe avec
l'axe des~$x_1$. La structure de contact~$\xi$ a alors pour \'equation de Pfaff
$$ \cos \theta \, dx_1 + \sin \theta \, dx_2 = 0, \qquad
   (x_1,x_2,\theta) \in \R^2 \times \R/2\pi\Z . $$
Le groupe $\DD (S,\partial S)$ est connexe (contractile d'apr\`es \cite{Sm})
et il faut donc montrer qu'il en est de m\^eme du groupe $\DD (V,\partial V ;
\xi)$. Pour cela, on note $D \subset V$ le disque m\'eridien $\{ \theta = 0 \}$
et on consid\`ere la fibration
$$ \DD (V,\partial V ;\xi) \longrightarrow \PP (D,V ;\xi) \,. $$
La fibre s'identifie au groupe $\DD (\D^3,\partial \D^3 ;\xi)$ qui est connexe
d'apr\`es le th\'eor\`eme~17. Par suite, il suffit de voir que l'espace $\PP
(D,V ;\xi)$ est connexe. Malheureusement, $D$ n'est pas convexe : le bord de~$V$
ne porte aucune courbe legendrienne isotope \`a $\partial D$. On recourt donc
\`a la proposition~8 ; le bord de~$V$ est en effet convexe et son feuilletage
caract\'eristique est scind\'e par la bi-courbe
$$ \Delta = \bigl\{ (\cos\varphi, \sin\varphi, \theta) \in \partial \D^2 \times
   \R/2\pi\Z \mid \sin (\theta-\phi) = 1 \bigr\} . $$
Pour plus de commodit\'e, on change les variables $x_1, x_2$ en
\begin{align*}
   x &= \cos\theta \, x_1 + \sin\theta \, x_2, \\
   y &= - \sin\theta \, x_1 + \cos\theta \, x_2 \,.
\end{align*}
On pose ensuite
\begin{align*}
   \alpha_0 &= dx - y\, d\theta = \cos\theta \, dx_1 + \sin\theta \, dx_2,  \\
   \alpha_1 &= (1-y^2) \, dx + xy \, dy - y \, d\theta, \\ \text{\llap{et\quad}}
   \alpha_s &= (1-s) \alpha_0 + s \alpha_1, \quad s \in [0,1] \,.
\end{align*}
Les structures de contact $\xi_s = \ker \alpha_s$ forment alors un chemin dans
$\SC(V,\Delta)$ issu de $\xi_0 = \xi$. Elles sont ainsi toutes tendues. De plus,
le disque $D$ est convexe dans $(V,\xi_1)$ et son feuilletage caract\'eristique
(figure~3) est scind\'e par la courbe $\Gamma = \{ y=0 \} \subset D$. Il reste
\`a prouver que $\PP(D,V ;\xi_1)$ est connexe. Compte tenu du lemme~7, de la
proposition~10 et de la connexit\'e de $\PP(D,V)$ \cite{Wa}, il suffit de montrer
que, si un plongement $\psi \in \PP (D,V ;\Gamma)$ a une image disjointe de~$D$
sauf au bord --~$D \cup \psi(D)$ borde alors automatiquement une ``boule
pinc\'ee''~--, il est isotope \`a l'inclusion parmi les plongements adapt\'es
\`a~$\Gamma$. Cela d\'ecoule de la classification des structures de contact
tendues sur la boule~\cite{El} : toute structure de contact tendue sur $\D^2
\times_\partial [0,1]$ pour laquelle les deux disques $\D^2 \times \{0\}$ et
$\D^2 \times \{1\}$ sont convexes est isotope, relativement au bord, \`a une
structure de contact pour laquelle tous les disques $\D^2 \times \{t\}$, $t \in
[0,1]$, sont convexes.

\begin{figure}[htbp]
\centering
\includegraphics[width=.5\textwidth]{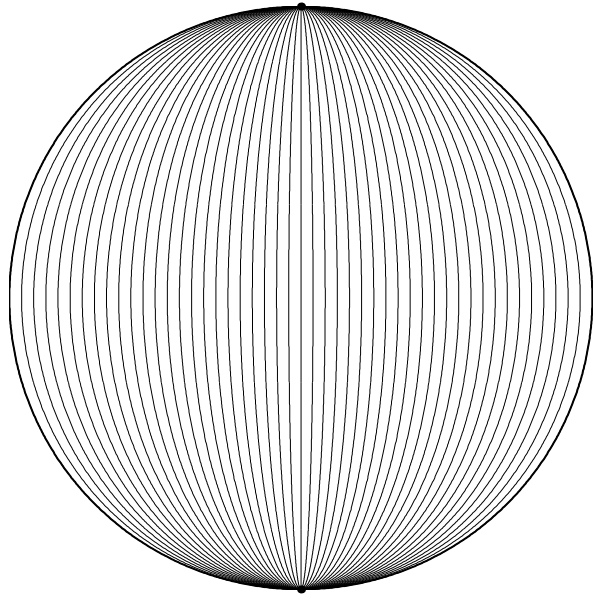}
\caption{Feuilletage caract\'eristique du disque m\'eridien~$D$ dans
$(V,\xi_1)$.}
\end{figure}

\medskip

Avant de d\'emontrer le th\'eor\`eme~1 pour les autres surfaces, on fait deux
observations g\'en\'erales.

\begin{lemme} 
Soit $C \subset S$ une  courbe simple, ferm\'ee ou propre. L'image inverse $F =
\pi^{-1}(C) \subset V$ est une surface convexe.
\end{lemme}

Dans la suite, on note $\Gamma \subset F$ la multi-courbe associ\'ee \`a un
voisinage homog\`ene quelconque de~$F$. La courbe $\Gamma$ a deux composantes
connexes et chacune d'elles se projette hom\'eomorphiquement sur~$C$.

\begin{proof}
Il existe sur $S$ un champ de vecteurs $\Ol X$ transversal \`a~$C$ dont le flot
plonge $C \times \R$ dans~$S$ en envoyant $\partial C \times \R$ dans $\partial
S$. Ce champ $\Ol X$, comme tout champ de vecteurs sur~$S$, se rel\`eve en un
champ de vecteurs~$X$ sur~$V$ qui pr\'eserve~$\xi_S$. L'image de $F \times \R$
par le flot de~$X$ est un voisinage homog\`ene de~$F$.
\end{proof}

\begin{lemme} 
Soit $C \subset S$ une  courbe simple essentielle, ferm\'ee ou propre, soit $F =
\pi^{-1}(C)$ et soit $\psi \in \PP(F,V ;\Gamma)$. Si $\psi(F)$ est isotope \`a
$F$, comme surface lisse et relativement au bord, la courbe $\psi(\Gamma)$ est
isotope \`a~$\Gamma$ relativement au bord. De plus, si $\psi$ est isotope \`a
l'inclusion dans $\PP (F,V)$, il l'est aussi dans $\PP(F,V ;\Gamma)$.
\end{lemme}

\begin{proof}
Soit $\Wt S$ le rev\^etement de~$S$ associ\'e \`a la courbe~$C$ et $\Wt C$ un
rel\`evement compact de~$C$ dans~$\Wt S$ ; si $C$~est un arc, $\Wt S$ est le
rev\^etement universel de~$S$. Par naturalit\'e, le fibr\'e $\wt\pi \from \Wt V
\to \Wt S$ des \'el\'ements de contact au-dessus de~$\Wt S$ n'est autre que le
rappel du fibr\'e $\pi \from V \to S$ au-dessus de~$\Wt S$. On note $\wt\xi$ la
structure de contact canonique de~$\Wt V$ et on pose $\Wt F = \wt\pi^{-1} (\Wt
C)$.

Toute isotopie $F_s$, $s \in [0,1]$, entre $F_0=F$ et $F_1=\psi(F)$ se rel\`eve
en une isotopie $\Wt F_s$ entre $\Wt F_0 = \Wt F$ et une surface $\Wt F_1$ qui
est convexe dans $(\Wt V,\wt\xi)$. De plus, on peut trouver dans $\Wt S$ un
domaine produit
$$ R = \Wt C \times \R \supset \Wt C = \Wt C \times \{0\}, \quad\text{avec}\quad
   \partial R = \partial \Wt C \times \R \subset \partial \Wt S, $$
qui renferme toutes les projections des surfaces $\Wt F_s$, $s \in [0,1]$. Le
domaine $U = \wt\pi^{-1}(R)$ est un voisinage homog\`ene de $\Wt F$ qui contient
$\Wt F_1$. L'in\'egalit\'e de Bennequin semi-locale [Gi3, proposition 4.10]
montre alors que les courbes qui scindent $\wt\xi \Wt F_1$ sont isotopes \`a
celles qui scindent $\wt\xi \Wt F$. Par suite, $\psi(\Gamma)$ est isotope \`a~$
\Gamma$.

On suppose maintenant que $\psi$ est isotope \`a l'inclusion dans $\PP(F,V)$.
Pour voir qu'il l'est aussi dans $\PP(F,V ;\Gamma)$, il suffit,  d'apr\`es la
proposition~10, de regarder le cas o\`u $\psi(F)$ est disjoint de~$F$ (sauf au
bord) et borde avec~$F$ un produit pinc\'e $W \simeq F \times_\partial [0,1]$.
Dans ce cas, les arguments ci-dessus montrent que $(W,\xi)$ se plonge dans le
voisinage homog\`ene $(U,\wt\xi)$. Il r\'esulte alors de la classification des
structures de contact sur le tore \'epais (si $C$ est une courbe ferm\'ee) et
sur le tore plein (si $C$ est un arc) que la structure de contact induite par
$\xi$ sur $F \times_\partial [0,1] \simeq W$ est isotope, relativement au bord,
\`a une structure de contact pour laquelle toutes les surfaces $F \times \{t\}$,
$t \in [0,1]$, sont convexes.
\end{proof}

On reprend maintenant la d\'emonstration du th\'eor\`eme~1 et on suppose que la
surface $S$ n'est ni une sph\`ere, ni un tore. On proc\`ede par r\'ecurrence sur
l'entier $n(S) = -2\Chi(S) - \beta(S)$ o\`u $\Chi(S)$ est la caract\'eristique
d'Euler de~$S$ et $\beta(S)$ le nombre de composantes connexes de $\partial S$.
On prend dans $S$ une courbe simple non s\'eparante~$C$, ferm\'ee ou propre, et
on note $S'$ la surface obtenue en d\'ecoupant~$S$ le long de~$C$ ; ainsi, $n(S')
< n(S)$.

Soit $\phi_0 \in \DD (V,\partial V ;\xi)$. Un th\'eor\`eme de F.~Waldhausen
\cite{Wa} assure que $\phi_0$ est isotope, dans $\DD (V,\partial V)$,  \`a un
diff\'eomorphisme fibr\'e $\phi_1$. On note $\ol\phi_1$ le diff\'eomorphisme de
$S$ sur lequel $\phi_1$ se projette et on pose $\phi = D(\ol\phi_1)^{-1} \circ
\phi_0$. En vertu des lemmes 18, 19 et~7, $\phi$ est isotope, dans $\DD (V,
\partial V ;\xi)$ \`a une transformation de contact~$\phi'$ qui induit
l'inclusion sur la surface $F = \pi^{-1}(C)$. En outre, la vari\'et\'e de
contact $(V',\xi')$ obtenue en d\'ecoupant $(V,\xi)$ le long de~$F$ s'identifie
\`a la vari\'et\'e des \'el\'ements de contact au-dessus de~$S'$. L'hypoth\`ese
de r\'ecurrence permet d\`es lors de conclure.

\subsection*
{D\'emonstration du th\'eor\`eme~1 : cas du groupe $\DD (V ;\xi)$}

Chaque diff\'eomorphisme $\phi \in \DD (V ;\xi)$ pr\'eserve le feuilletage
caract\'eristique de $\partial V$, donc les fibres de la projection $\partial V
\to \partial S$ qui sont legendriennes. Comme toute isotopie de $\partial V$ qui
laisse invariant le feuilletage caract\'eristique s'\'etend en une isotopie de
contact, on peut d\'eformer~$\phi$ dans $\DD(V ;\xi)$ en une transformation qui
co\"{\i}ncide, sur $\partial V$, avec la diff\'erentielle d'un diff\'eomorphisme
de~$S$. On applique alors le r\'esultat sur $\DD (V,\partial V ;\xi)$ pour
conclure.

\subsection*
{D\'emonstration du th\'eor\`eme~1 : cas du groupe $\DD (\Int V ;\xi)$}

Soit $\phi \in \DD (\Int V ;\xi)$. Les r\'esultats de F.~Waldhausen~\cite{Wa}
joints \`a l'in\'egalit\'e de Bennequin semi-locale~\cite{Gi3} assurent, comme
dans le lemme~19, que $\phi$ est isotope dans $\DD (\Int V)$ \`a la
diff\'erentielle d'un diff\'eomorphisme de~$S$.  Quitte \`a composer $\phi$ avec
l'inverse de celle-ci, on suppose que $\phi$ est isotope \`a l'identit\'e dans
$\DD (\Int V)$. On choisit dans~$S$ un voisinage collier du bord
$$ R = \partial S \times [0,\infty \cl[ \supset
   \partial S = \partial S \times \{0\} $$
et on pose $U = \pi^{-1}(R)$. Le champ de vecteurs $\partial_t$, $t \in \R$,
donn\'e sur~$R$ par la structure produit se rel\`eve en un champ de vecteurs de
contact qui fait de~$U$ un voisinage collier homog\`ene de $\partial V$. Pour
$t_0>0$ assez petit, l'image par $\phi$ du multi--tore
$$ F = \partial V \times \{t_0\} \subset U = \partial V \times [0,\infty\cl[ $$
est contenue dans~$U$ et les m\^emes arguments que dans la d\'emonstration du
lemme~19 et de la proposition~8 montrent que $\phi$ est isotope, dans
$\DD (\Int V ;\xi)$, \`a une transformation de contact~$\phi'$ qui co\"\i ncide
avec l'identit\'e sur un voisinage de~$F$. On choisit alors sur $S$ un champ de
vecteurs $\Ol X$ ayant les propri\'et\'es suivantes :

\begin{itemize}
\item
$\Ol X = 0$ hors de $\partial S \times [0,t_0+\eps] \subset R$ o\`u $\eps>0$ est
pris assez petit pour que $\phi'$ induise l'identit\'e sur $\partial V \times
[t_0,t_0+\eps]$ ;
\item
$\Ol X = f(t) \partial_t$ sur $\partial S \times [0,t_0+\eps]$ o\`u $f$ est une
fonction strictement n\'egative sur l'int\'erieur de $[t_0,t_0+\eps]$ et nulle
au bord.
\end{itemize}
On note $X$ le champ de vecteurs de contact sur~$V$ qui rel\`eve $\Ol X$ et
$\tau_s$, $s \in \R$, le flot de~$X$. Les transformations de contact
$$ \phi_s = \tau_s \circ \phi' \circ \tau_{-s}, \qquad s \in [0,\infty], $$
constituent un chemin dans $\DD (\Int V ;\xi)$ qui relie $\phi_0 = \phi'$ \`a
un diff\'eomorphisme \'egal \`a l'identit\'e pr\`es du bord. On applique alors
le r\'esultat sur $\DD (V,\partial V ;\xi)$ pour conclure.

\subsection*{D\'emonstration du th\'eor\`eme~3}

Soit $C$ une courbe ferm\'ee simple sur $S = \S^2$, soit $F$ le tore convexe
$\pi^{-1} (C) \subset (V,\xi)$ et soit $\Gamma$ une multi-courbe qui scinde $\xi
F$. Le tore~$F$ partage $V \simeq \SO_3$ en deux tores pleins $W_\pm$ et chaque
composante connexe de~$\Gamma$ rencon\-tre homologiquement une fois le m\'eridien
de~$W_\pm$. Soit maintenant $\phi \in \DD (V ;\xi)$. D'apr\`es le th\'eor\`eme
de J.~Cerf~\cite{Ce}, $\phi$ est isotope \`a l'identit\'e dans $\DD(V)$. Pour
voir qu'il l'est aussi dans $\DD (V ;\xi)$, il suffit, vu le th\'eor\`eme~1, de
montrer que $\phi \res F$ est isotope \`a l'inclusion dans $\PP(F,V ;\Gamma)$.
Compte tenu de la proposition~10, on peut supposer que $\phi(F)$ est disjoint
de~$F$ et borde avec~$F$ un produit $W \simeq F \times [0,1]$. Il d\'ecoule
alors de la classification des structures de contact sur le tore \'epais~\cite{Gi2}
que $\phi(F)$ est isotope \`a~$F$ parmi les tores convexes.

\subsection*{D\'emonstration du th\'eor\`eme~4}

Le fait que l'image du morphisme
$$ \pi_0 \DD (\T^3,\xi) \longrightarrow \SL_3(\Z) $$
soit exactement le groupe~$\Pi$ est d\'emontr\'e dans~\cite{EP}. L'injectivit\'e
du morphisme r\'esulte du lemme~19 comme le th\'eor\`eme~1.

\adresse{
Emmanuel \textsc{Giroux}\ls
Unit\'e de Math\'ematiques Pures et Appliqu\'ees\ls
\'Ecole Normale Sup\'erieure de Lyon\ls
46, all\'ee d'Italie\ls
69364 Lyon cedex 07, France
}

\end{document}